\documentclass[11pt]{article} 

\usepackage[utf8]{inputenc} 

\usepackage{geometry} 
\usepackage{graphicx} 

\usepackage{booktabs} 
\usepackage{array} 
\usepackage{paralist} 
\usepackage{verbatim} 
\usepackage{subfig} 

\usepackage{fancyhdr} 
\usepackage{sectsty}
\allsectionsfont{\sffamily\mdseries\upshape} 

\usepackage[nottoc,notlof,notlot]{tocbibind} 
\usepackage[titles,subfigure]{tocloft} 

\usepackage{amsmath, amsfonts, amsthm, amssymb, amscd}
\usepackage{indentfirst}

\newtheorem{theorem}{Theorem}

\title{A Note on Counting Dependency Trees}
\author{Zhujun Zhang}
         
\providecommand{\keywords}[1]{\textbf{\textit{Index terms---}} #1}

\begin{document}
\maketitle

\begin{abstract}
We apply symbolic method to deduce functional equation which generating function of counting sequence of dependency trees must satisfy. Then we use Lagrange inversion theorem to obtain concrete expression of the counting sequence. We apply the famous Stirling's approximation to get approximation of the counting sequence.  At last, we discuss the additive parameters of dependency trees.
\end{abstract}


\section{Introduction}
Dependency is a one-to-one correspondence: for every element (e.g. word or morph) in the sentence, there is exactly one node in the structure of that sentence that corresponds to that element. The result of this one-to-one correspondence is that dependency grammars are word (or morph) grammars. All that exist are the elements and the dependencies that connect the elements into a structure. The structure could be represented by dependency trees(cf. [2]). 

Counting dependency trees is a foundation of average-case analysis of algorithms on dependency grammars processing. Hu et al. found a recurrence counting formula of dependency trees in [1]. Marco Kuhlmann deemed a closed form of counting sequence of dependency trees in [3].

In this note, we apply symbolic method to deduce counting formula of dependency trees. And we use famous Stirling's approximation to obtain approximation of counting sequence.

\section{Preliminaries}
Dependency trees could be defined recursively as: a dependency tree is either (i) a single node tree, or (ii) a root node with several sub-dependency trees on left and right respectly.

Symbolic method could be considered as a set of intuitive combinatorial constructions that immediately translate to equations that the associated generating functions must satisfy. Symbolic method provides translation tools for a large number of combinatorial constructions on combinatorial classes. Here, we just introduce some elementary translation tools. Advanced symbolic method and extensive applications could be found in Chapter 5 of [4] and Part A of [5].

To specify combinatorial classes, we make use of neutral objects $\epsilon$ of size 0 and the neutral class $\mathcal{E}$ that contains a single neutral object. 
Then we introduce three simple operations on combinatorial classes.
Given two classes $\mathcal{A}$ and $\mathcal{B}$ of combinatorial objects, we can build new classes as follows:
$$
\begin{aligned}
	& \mathcal{A} + \mathcal{B} \ \text{is the class consisting of disjoint copies of the members of} \ \mathcal{A} \ \text{and} \ \mathcal{B}, \\
	& \mathcal{A} \times \mathcal{B} \ \text{is the class of ordered pairs of objects, one from} \ \mathcal{A} \ \text{and one from} \ \mathcal{B}, \\
	& \text{and} \ SEQ(\mathcal{A}) \ \text{is the class} \ \epsilon + \mathcal{A} + \mathcal{A} \times \mathcal{A} + \mathcal{A} \times \mathcal{A} \times \mathcal{A} + \dots \ \text{.}
\end{aligned}
$$

Following elementary symbolic method provides a simple correspondence between operations in combinatorial constructions and their associated generating functions.
\begin{theorem}
	\textbf{(Symbolic method)}. 
	Let $\mathcal{A}$ and $\mathcal{B}$ be unlabelled classes of combinatorial objects.
	If A(z) is the generating function that enumerates $\mathcal{A}$ 
	and B(z) is the generating function that enumerates $\mathcal{B}$, then
	$$
	\begin{aligned}
	A(z) + B(z) \  & \text{is the generating function that enumerates} \  \mathcal{A} + \mathcal{B} \\
	A(z)B(z) \  & \text{is the generating function that enumerates} \  \mathcal{A} \times \mathcal{B} \\
	\frac{1}{1 - A(z)} \  & \text{is the generating function that enumerates} \  SEQ(\mathcal{A}).
	\end{aligned}
	$$
\end{theorem}


Lagrange inversion theorem is of particular importance for tree enumeration. The theorem allows us to extract coefficients from generating functions that are implicitly defined through functional equations.

\begin{theorem}
\textbf{(Lagrange inversion theorem)}. Suppose that a generating function $A(z)=\sum_{k \geqslant 0} {a_k z^k} $ satisfies the functional equation 
$z=f(A(z))$ , where $f(z)$ satisfies $f(0)=0$ and $f'(0) \neq 0$. Then
$$
\begin{aligned}
a_n \equiv \left [ z^n \right ] A(z)  = \frac{1}{n} \left [ u^{n-1} \right ] {\left (\frac{u}{f(u)} \right )}^n.
\end{aligned}
$$
\end{theorem}

\section{Couting Dependency Trees}
Let $\mathcal{T}$ be combinatorial class of dependency trees. Definition of dependency trees implies the construction of corresponding combinatorial class
$$
\begin{aligned}
\mathcal{T} = SEQ(\mathcal{T}) \times \circ \times SEQ(\mathcal{T})
\end{aligned},
$$
where $\circ$ be denoted as root node of dependency trees. Let $T(z)$ be generating function of $\mathcal{T}$, and generating function of single node is $z$. Symbolic method could translate the construction to functional equation
$$
\begin{aligned}
T(z) = \frac{1}{1-T(z)} \cdot z \cdot \frac{1}{1-T(z)}
\end{aligned}.
$$

With the help of \textbf{WolframAlpha}, we can get exact expression of $T(z)$, 
$$
T(z)= \frac{1}{3}\left( \frac{\sqrt[3]{3\sqrt{3} \sqrt{27z^2-4z} +27z -2}  }{\sqrt[3]{2}} + \frac{\sqrt[3]{2}}{\sqrt[3]{3\sqrt{3} \sqrt{27z^2-4z} +27z -2}  } +2\right) .
$$
However, it is difficult to extract coefficient of generating function from the expression above. Fortunately, Lagrange inversion theorem could be applied. Through elementary algebra, we obtain $T(z)(1-T(z))^2=z$. Let $f(u)=u(1-u)^2$, we have 
$$
\begin{aligned}
t_n\equiv\left [ z^n \right ] T(z) & = \frac{1}{n} \left [ u^{n-1} \right ] \frac{1}{(1-u)^{2n}} \\
& = \frac{1}{n} \left [ u^{n-1} \right ] \sum_{k \geqslant 0} {\binom{k+2n-1}{k} u^k} \\
& = \frac{1}{n} \binom{n-1+2n-1}{n-1} \\
& = \frac{1}{n} \binom{3n-2}{n-1} .
\end{aligned}
$$

The sequence is \textbf{OEIS A006013}(cf. [3]).
Then, with Stirling's approximation we could obtain approximation of $t_n$ 
$$
\begin{aligned}
t_n & = \frac {n(2n)3n!} {n (3n)(3n-1)n!2n!} \\
& \sim \frac{2}{9n} \frac{ \sqrt{6 \pi n} \left ( 3n / e \right ) ^{3n} }
{\sqrt{2 \pi n} \left ( n / e \right ) ^{n} \sqrt{4 \pi n} \left ( 2n / e \right ) ^{2n} } \\
& =  \frac {1} { \sqrt {27\pi} n^{3/2}} \left( \frac {27} {4} \right) ^n .
\end{aligned}
$$
Thus, we can conjecture $z_0=4/27$ is the dominant singularity of generating function $T(z)$.

\section{Additive Parameters of Dependency Trees}
We could defne an additive parameter to be any parameter whose cost function satisfes the linear recursive schema
$$
c(t)=e(t)+ \sum_{r \in t_s} c(r) ,
$$
where the sum is over all the subtrees of the root of $t$. The function $e$ is called toll function. Let $C(z)$ be cumulative generating function of an additive tree parameter $c(t)$ for the dependency trees, and let $E(z)$ be cumulative generating function for the associated toll function $e(t)$. We could deduce the relation between these two functions
$$
\begin{aligned}
C(z) &\equiv \sum_{t \in \mathcal{T}} {c(t)z^{\left| t \right| }} \\
&= \sum_{t \in \mathcal{T}} {e(t)z^{\left| t \right| }} + \sum_{t \in \mathcal{T}} z^{\left| t \right| } \sum_{r \in t_s} c(r) \\
&= E(z) + \sum_{k \geqslant 0} (k+1) \sum_{t_1 \in \mathcal{T}} \cdots \sum_{t_k \in \mathcal{T}} \left( c(t_1) + \cdots + c(t_k) \right) z^{\left| t_1 \right| + \cdots + \left| t_k \right| +1} \\
&= E(z) + z\sum_{k \geqslant 0} (k+1) k C(z) T^{k-1}(z) \\
&= E(z) + \frac{2zC(z)}{ \left( 1- T(z) \right) ^3} .
\end{aligned}
$$
Thus, we obtain the equation of cumulative generating functions.
$$
C(z)= \frac{E(z)}{1-2z/\left( 1- T(z) \right) ^3} =E(z) \left( \frac{1-T(z)}{1-3T(z)} \right) .
$$


\begin{thebibliography}{12}
\bibitem {1} Fengguo Hu, Wei Huang, Haitao Liu, Enumeration of Dependency Structural Trees, Computer Engineering and Applications(in Chinese), 2009, Vol. 45(32), pp.22-24.
\bibitem {2} "Dependency grammar" of Wikipedia, https://en.wikipedia.org/wiki/Dependency\_grammar.
\bibitem {3} Neil J. A. Sloane, "A006013" of 
The On-Line Encyclopedia of Integer Sequences(OEIS), http://oeis.org/a006013.
\bibitem {4} Robert Sedgewick, Philippe Flajolet, An Introduction to the Analysis of Algorithms(2nd edition), Addison-Wesley Professional, 2013.
\bibitem {5} Philippe Flajolet, Robert Sedgewick, Analytic Combinatorics, Cambridge University Press, 2009.
\bibitem {6}  Ronald L. Graham, Donald E. Knuth, Oren Patashnik, Concrete Mathematics(2nd edition), Addison-Wesley Professional, 1994.
\end{thebibliography}
\end{document}